\documentclass[11pt]{article}
\usepackage{graphicx}
\usepackage{amsmath, amsfonts}

\title{On the (ab)use of statistics in the legal case against the nurse Lucia de 
B.}
\author{Ronald Meester\footnote{Vrije Universiteit Amsterdam}, 
Marieke Collins\footnote{Universiteit Utrecht}, Richard Gill\footnote{Universiteit 
Utrecht}, Michiel van Lambalgen\footnote{Universiteit van Amsterdam}}

\begin{document}
\maketitle

\begin{abstract}
We discuss the statistics involved in the legal case of the nurse Lucia de B.\ in 
The Netherlands, 2003-2004. Lucia de B. witnessed an unusually high number of incidents 
during her shifts, and the question arose as to whether 
this could be attributed to chance. We discuss and criticise 
the statistical analysis of Henk Elffers, a statistician 
who was asked by the prosecutor to write a statistical 
report on the issue. We discuss several other possibilities 
for statistical analysis. Our main point is that several 
statistical models exist, leading to very different predictions, or 
perhaps different answers to different questions. There is no 
such thing as a `best' statistical analysis.
\end{abstract}

\section{Introduction; the case} In The Hague (The Netherlands), 
on March 24, 2003 the nurse Lucia de B.\ (hereafter called 
either `Lucia' or `the suspect') was sentenced to life 
imprisonment for allegedly killing or attempting to kill a number of patients in 
two hospitals where she had worked in the recent past: the {\em Juliana 
Kinderziekenhuis} (JKZ) and the {\em Rode Kruis Ziekenhuis} (RKZ). 
At the RKZ she worked in two different wards, numbered 41 and 42 respectively.
At the JKZ, an unusually high proportion of incidents occurred 
during her shifts,\footnote{The precise technical definition of
`incident' is not important here; suffice it to say that an incident refers to 
the  necessity for reanimation, regardless of the outcome of the reanimation.} 
and the question arose as to whether Lucia's 
presence at so many `incidents' could  
have been merely coincidental.
 
A statistical analysis was given by statistician Henk Elffers, who 
had looked into the matter at the request of the public prosecutor. In broadest 
terms, his conclusion was this: assuming only (as he says) that
\begin{enumerate}
\item the probability that the suspect experiences an incident during 
a shift is the same as the corresponding probability for any other nurse,
\item the occurrences of incidents are independent for different shifts,
\end{enumerate}
then the probability that the suspect has experienced as 
many incidents as she in 
fact has, is less than $1$ in $342$ million.  According to Elffers, this 
probability is so small that standard statistical methodology sanctions 
rejection of the null-hypothesis of chance. He did take care to note 
that in itself this does not mean the suspect is guilty.

Why do we write this article? Two of us (MvL and RM) became involved in the case 
as expert witnesses of the defence. We studied the method and conclusion of Elffers,
and came to the conclusion that his numbers did not mean very much, if anything at all.
Elffers (and the court, for that matter) completely overlooked the subjective
element in the choice of a probabilistic model, and therefore the possibility 
of there being several models with very different predictions, or perhaps different
answers to different questions! 

The question as to how to use statistics in a case 
like this is not a question with a well-defined answer. 
Borrowing a phrase of Thomas Kuhn, we deal here with a {\em problem} 
rather than with a {\em puzzle}.
There are many ways of doing statistics. One can argue whether
to use a (subjective) Bayesian approach, or a classical frequentist
approach. There is even a school called the likelihood approach
which says that you should compute and report likelihood ratios,
full stop. Within each school there can be many solutions to what
appears to be the same problem. Moreover there is the question of
the range of the model.

Hence, many different approaches are possible, 
using very different models, and with many
different levels of sophistication. One can 
choose a very simple model, as Elffers did,
giving precise results, albeit of limited 
relevance. One can also choose a much
broader perspective, like a Bayesian point of 
view, which involves much more data, but
whose conclusions are much less precise. 
There simply is no unique best way of dealing with
the problem, and in this paper we want 
to elaborate on this point significantly. In court,
the judges continued to ask us: ``So if you 
reject Elffers' numbers, why don't you give us
better numbers'', implicitly assuming that 
there exist something like best numbers. 
One of the points of the present article is 
to argue against this.

This article is structured as follows. 
We will first present the relevant data and the 
statistical methodology used by Elffers. We discuss and 
criticise this methodology on several levels: not
only do we offer a critical discussion of his overall 
approach, but we also think that within his
paradigm, Elffers made several important mistakes. 
We also briefly discuss the way the court 
interpreted Elffers' report. Then
we show how the method of Elffers could have 
been used in a way we believe is correct
within his chosen paradigm, leading to a very 
different conclusion. After that, we discuss 
a Bayesian point of view, as advocated 
by the Dutch econometrician De Vos, and then we 
move on to the so called epidemiological models,
inspired by recent work of Lucy and Aitken. 
In the final section we try to formulate some conclusions.

\section{The data and Elffers' method} 

Elffers tried to base his model entirely on data pertaining to 
shifts of Lucia and the other nurses, and the incidents 
occurring in those shifts. The data on shifts 
and incidents for the period which was singled out 
in Elffers' report are given in the following table:
 
\medskip 

\medskip 

\medskip

\begin{tabular}{|l|c|c|c|}
\hline
hospital name (and ward number)& JKZ & RKZ-41 & RKZ-42\\ \hline
total number of shifts & $1029$ & $336$ & $339$ \\ \hline
Lucia's number of shifts & $142$ & $1$& $58$\\ \hline \hline
total number of incidents & $8$ & $5$ & $14$ \\ \hline
number of incidents during Lucia's shifts & $8$ & $1$ & $5$ \\ \hline 
\end{tabular}

\medskip \medskip
Later it was discovered that Lucia actually had done $3$ shifts in RKZ-41 
instead of just $1$, and in our own computations later in this article, 
we will use this correct number. 

When trying to put the situation sketched into a statistical model, 
ones first choice might be to build a model on the basis 
of epidemiological data concerning the probability of incidents during 
various types of shifts; this would allow one to calculate the probability 
that the suspect would be present accidentally at as many incidents 
as she in fact witnessed. 
 
However, the trouble with this approach is that for the most part 
the requisite data are lacking. And even if the data were available, their
use would be a subject of debate between prosecutor and defence; 
see Section \ref{ait}.

Because of this, Elffers tried to set up a model which uses only the
shift data given above. This he achieved by {\em conditioning} on part 
of the data. He assumed that

\begin{enumerate}
\item there is a fixed probability $p$ for the occurrence 
of an incident during a shift (hence $p$ does not depend on whether the
shift is a day or a night shift, etc.),
\item incidents occur independently of each other.
\end{enumerate}

It is now straightforward to compute the {\em conditional} probability 
of the event that (at the JKZ, say) all incidents occur during Lucia's shifts,
{\em given} the total number of incidents and the total number of shifts 
\textit{in the period under study}.
Indeed, if the total number of shifts is $n$, and Lucia had $r$ shifts,
then the conditional probability that Lucia witnessed $x$ incidents given that
$k$ incidents occurred, is
\begin{equation}
\label{hypergeometric}
\frac
{\displaystyle 
{r \choose x} p^x (1-p)^{r-x} ~ { n-r \choose k-x} p^{k-x}(1-p)^{n-r-k+x}
}
{\displaystyle 
{n \choose k} p^k(1-p)^{n-k}
}
~=~
\frac
{\displaystyle 
{r \choose x}{n-r \choose k-x} 
}
{\displaystyle {n \choose k}
}.
\end{equation}
Note that this quantity does not depend on the unknown parameter $p$. This
distribution is known as the {\em hypergeometric} distribution. With this 
formula, one can easily compute the (conditional) probability that the suspect 
witnessed at least the number of incidents as she actually has, for each ward. 

However, according to Elffers, this computation is not completely 
fair to the suspect.
Indeed, the computation is being done precisely because there
were so many incidents during her shifts at the JKZ. It would, 
therefore, be more reasonable (according to Elffers) not
to compute the probability that Lucia has witnessed so many 
incidents, but instead the probability that {\em some} nurse 
witnessed so many incidents. At the JKZ, there were
$27$ nurses taking care of the shifts and therefore, presumably
to get an upper bound to this probability, Elffers multiplies his
outcome by $27$; he calls this the {\em post hoc correction}. 
According to Elffers, this correction only
needs to be done at the JKZ; at the RKZ this is no longer
necessary since the suspect was already identified as being 
suspect on the basis of the JKZ data.

Elffers arrives at his final figure (the aforementioned $1$ in $342$ million) 
by multiplying the outcomes for the three wards (with post hoc correction 
at the JKZ, but without this correction at the RKZ).

\section{Discussion of Elffers' method} 

There are a number of problems and points of concern 
with the method of Elffers. In the following, we
list some of these. 

\subsection{Conditioning on part of the data}
As we remarked already, conditioning on the number of incidents has a big advantage,
namely that under the hypothesis of chance, the unknown parameter $p$ cancels in
the computations. It is the very conditioning that makes computations possible
in Elffers' model. 

The idea of conditioning at inference time on quantities that
were not naturally fixed at data sampling time has some history. It seems that
Fisher first proposed this idea for exact inference on a $2 \times 2$ contingency
table \cite{fisher}. In \cite{mpatel}, some justification is offered for this 
technique. Conditioning is reasonable, according to Mehta and Patel, if ``the
margins [...] are representative of nuisance parameters whose values do not
provide any information about the null-hypothesis of interest." (They mean:
information about the truth or falsity of the hypothesis).
When we discuss loss of information by conditioning, the generally accepted 
attitude is that the loss is usually not worth fussing about. 
However, we deal here with a legal case, and the fact 
that {\it usually} the information loss might be not worth fussing about, is
not enough to dismiss this issue in an individual legal case. 
It is also usually the case that DNA material on the body of a crime victim comes
from the criminal, but in each individual case, this needs to be confirmed.

In the current case, it is not clear at all that the number 
of incidents does not provide information.
We do not know for sure that \textit{any} murders were committed;
all other evidence was circumstantial.
If the total number of incidents happened to be smaller 
than normally would be expected, it becomes less plausible 
that many attempted murders have taken place.

\subsection{Using data twice: the post hoc correction}
One of the problems with this approach is the fact that the data of the JKZ 
is used twice. First to identify the suspect and indeed, 
to suspect that a crime has occurred, and after that again in the computations
of Elffers' probabilities. This procedure should raise eyebrows amongst statisticians:
it is one of these problems that seem to arise all over the place: one sets up
an hypothesis on the basis of certain data, and after that one uses the same
data to test this hypothesis. It is clear that this raises problems, and it is equally
clear that Elffers' method shares this problem. In a way, Elffers seems to be aware
of this. After all, his post hoc correction was introduced for exactly this reason.

However, this post hoc correction is a striking instance of an 
unacknowledged subjective choice employed by Elffers. 
To see this, note that Elffers restricts the statistical analysis 
to the wards at which the suspect worked. Why?
The question of the prosecutor, whether Lucia's number of incidents
can be put down to chance, has to do with the question whether or
not Lucia killed or attempted to kill some of her patients. 
The word `ward' is not mentioned in these questions, nor is there 
any mention of the other nurses who worked there. 
It was the choice of Elffers himself to consider the level of 
wards. We do not claim that this decision was wrong; 
there are arguments to defend it, the most important one
probably being the simplicity of the resulting model. 
But one can also envision a statistical analysis of {\em all} 
wards in, say, The Netherlands, perhaps with different
probabilities for incidents in different wards. We might now condition 
on the number of incidents in each ward. Again, the number of incidents 
of the suspect has the same hypergeometric distribution as before. 
But the necessary post hoc correction in this hypothetical statistical analysis 
would logically take account of {\em all} nurses in The Netherlands, 
even though the computations concerning the suspect might still only 
depend on the data of her own ward. Multiplication by the number of nurses in the
ward of Lucia does not necessarily follow from the fact that we only use data
from her own ward; the level of the post hoc correction is arbitrary.

An analogy might clarify this point.
Consider a lottery with tickets numbered $1$ to $1,000,000$. The jackpot 
falls on the ticket with number $223,478$, and the ticket has been bought 
by John Smith. John Smith lives in the Da Costastraat in the city of Leiden. 
Given these facts we may compute the chance that John Smith wins the jackpot; 
a simple and uncontroversial model shows that this probability will be
extremely small. Do we conclude form this that the lottery was not fair, 
since an event with very small probability has happened? Of course not. We can
also compute the probability that someone in the Da Costastraat wins the jackpot, 
but it should be clear that the choice of the Da Costastraat as reference 
point is completely arbitrary. We might similarly compute the probability 
that someone in Leiden wins the jackpot, or someone living in Zuid-Holland 
(the state in which Leiden is situated). With these data-dependent hypotheses 
there simply is no uniquely defined scale of the model at which the problem must be
studied. 

The analogy with the case of Lucia will be clear: the winner of
the jackpot represents the suspect being present at $8$ out of $8$ 
incidents, the street represents the ward. Elffers restricts his model 
to the ward in which something unusual has happened. With perhaps equal justification, 
another statistician might have considered the entire JKZ (Leiden, in the 
analogy) instead of the ward as basis for her computations -- with vastly higher 
probability for the relevant event to happen somewhere. Still another statistician might 
have taken the Netherlands as the basis for the computation, which yields again 
a higher probability. The important point to note is that \textit{subjective
choices are unavoidable here}; and it is rather doubtful whether a court's 
judgement should be based on such choices. If one wants to avoid this kind of
subjective choice, one should adopt an approach where the data is not used twice. 
In the next section we discuss such an approach.

Even if we agree with the level (wards) of the posthoc correction,
the way it is done needs to be motivated. Elffers' motivation 
is to compute (though presumably, he means to bound)
the probability that some nurse among the $27$ at JKZ 
would experience as many incidents as Lucia. 
A glance at the numbers shows that the $27$ nurses must have had
very varied numbers of shifts. The chance that any particular nurse
would experience at least $14$ incidents will depend on her total
number of shifts, and appears hardly relevant. More relevant would
perhaps be each nurse's incident rate or risk 
(average number of incidents per shift), but we are not given the
numbers of shifts of the other nurses.

\emph{If} we suspected a priori that murders had taken place at the JKZ, 
and wanted to investigate whether they were associated with one of the
nurses, then \textit{before seeing the data} a statistician \emph{might} 
reasonably adopt the following standard (Bonferoni)
procedure for so-called \textit{multiple comparisons}:
compute for each nurse separately, the probability of their witnessing
at least the number of incidents which they did witness, under the hypothesis
of randomness. Multiply the smallest of these probabilities by the number
of nurses. The result would be a legitimate $p$-value (the meaning of 
`$p$-value' is discussed in the next section). Thus Elffers'
post hoc correction \emph{could} have been appropriate under a
rather different unfolding of the events. But this does not
justify his correction in the present circumstances. 

Other commentators have derived Elffers' post hoc correction in a
Bayesian approach where it is assumed that there have been murders
by a nurse, and that each nurse has an equal probability of being the 
murderer. Again, this picture simply does not apply to the actual 
circumstances of the case.

\subsection{Multiplication is not allowed}
\label{mult}

Elffers multiplies the three probabilities from the three wards. 
The multiplication means that he is assuming that under his 
null-hypothesis, incidents occur completely randomly in each of the
three wards (as far as the allocation of shifts to nurses is 
concerned), independently over the wards, 
but with possibly different rates in each ward. 
If one accepts his earlier null-hypothesis as an
interesting hypothesis to investigate, then this new hypothesis 
could also be of interest.

What is the meaning of the probability which Elffers finds? 
It is the probability, under this null-hypothesis of randomness, 
and conditional on the total number of incidents 
in each ward, that a nurse with as many shifts as Lucia in each ward
separately, would experience as many (or more) incidents than she did, 
in all wards simultaneously. 
Is the fact that this probability is very small, good reason 
to discredit the null hypothesis?

First we should understand the rationale of Elffers' method when applied 
to one ward. He is interested to see if a certain null-hypothesis is 
tenable (whether his null-hypothesis is relevant to the 
case at hand, is discussed in the next section). 
He chooses in advance for whatever reason he likes, a 
statistic (a function of the data) such that large values of that 
statistic would tend to occur more easily if there actually is a, 
for him, interesting deviation from the null-hypothesis. 
Since his null-hypothesis completely fixes the distribution of his
chosen statistic, he can compute the probability that the actually 
observed value could be equalled or exceeded under that hypothesis. 
The resulting probability is called the $p$-value of the statistical test. 
If null-hypothesis and statistic are specified in advance of 
inspecting the data, then it becomes hard to retain belief in the 
null-hypothesis if the $p$-value is very small. Elffers in fact 
follows the following procedure: he has selected (arbitrarily) 
a rather small threshold, $p=0.001$. When a $p$-value is smaller 
than $0.001$ he will declare that the null-hypothesis is not true.
Following this procedure, and in those cases when actually the null-hypothesis 
was true, he will make a false declaration once in a thousand times.

If the null-hypothesis corresponds to a person being innocent of having committed a
crime, then his procedure carries the guarantee that not more than one in a
thousand innocent persons are falsely found guilty. (Presumably, society does
accept some small rate of false convictions, since absolute certainty about guilt
or innocence is an impossibility. But perhaps one in a thousand is a
bit too large a risk to take).

Now we return to Elffers' multiplication of three $p$-values, one for each ward. 
Does this result in a new $p$-value?

An easy argument shows that the answer is \textit{no}. Suppose
there are $100$ wards and the null-hypothesis is true (including the independence 
over the wards).
A nurse with the same number of shifts as Lucia in each ward 
has approximately a probability
of a half to have as many incidents as Lucia, in each ward separately. 
Multiplying, the probability 
that she `beats' Lucia in all wards is approximately 1 in 2 to the power one hundred, 
or approximately one in a million million million million million. Yet we are assuming 
the complete randomness of incidents within each ward! 
Clearly we have to somehow discount
the number of multiplications we are doing.

Is there something else that Elffers could have done, to combine the 
results of the three wards?
Yes; and in fact, classical statistics offers many choices.
For instance he could have compared the total number of incidents
of Lucia over the three wards, to the probability of exceeding that number, 
given the totals per ward
and the numbers of shifts, when in each ward separately 
incidents are assigned uniformly at
random over all the shifts. In the language of statistical hypothesis testing, 
he should have chosen
a single test-statistic based on the combined data for his combined null-hypothesis, 
and computed its $p$-value under that null-hypothesis. Perhaps it would be reasonable
to weight the different wards in some way. 
Each choice gives a different test-statistic and a different
result. The choice should be made in advance of looking at the data, 
and should be designed
to react to the kind of deviation from the null-hypothesis which it is most important to
detect. Such a choice can be scientifically motivated but it is in the last analysis 
subjective.

An easy way to combine (under the null-hypothesis) 
independent $p$-values is a method due to Fisher 
(which can be found in his book 
\cite{fisher}): 
multiply (as Elffers did) the three $p$-values
for the separate tests (denoted by $p_1, p_2$ and $p_3$), 
and compare this with the probability distribution of the product 
of the same number of uniform random numbers between 0 and 1. 
A standard argument from probability
theory reduces this to a comparison of  $- 2 \sum_i \log p_i$   
with a chi-squared distribution 
with $2n$ degrees of freedom, where $n=3$. 
What is in favour of this method is its simplicity. 
Choosing this one is just as much a subjective choice as any other.
The final $p$-value will be much larger than that reported by
Elffers.

\subsection{The Quine-Duhem problem}
\label{qd}

In fact, talk of
`the rejection of the null-hypothesis' is somewhat imprecise. It
was observed by the philosopher Quine, and before him by the
historian of science Duhem, that the falsificationist picture of 
an hypothesis $H$ logically implying a prediction $P$, which when 
falsified must lead to the abandonment of $H$, is too simplistic. 
 
Consider the following example. Suppose our 
thermodynamic theory implies that water 
boils at 100C at sea level; and suppose furthermore that our observations
show water to boil at 120C. Does this mean thermodynamics is false?
Not necessarily, because there might be something wrong with the
thermometer used. That is, the logical structure of the prediction is rather 
\begin{quote} `Thermodynamics + Properties of thermometer' imply 
`water boils at 100C at sea level'.
\end{quote}
More formally, a prediction $P$ from an hypothesis $H$ always has 
the form $H\&A_1\&\cdots\& A_n \Rightarrow P$, where the $A_1\&\cdots\& A_n$ 
are the auxiliary hypotheses. If we find that $P$ is false, we can 
conclude only not-$(H\&A_1\&\cdots\& A_n)$ from which something can be concluded 
about $H$ only if we have independent corroboration of the $A_1\&\cdots\& A_n $.
The same phenomenon occurs in statistics, and in particular in this case.

In order to be able to make calculations, Elffers
in fact explicitly makes auxiliary assumptions far beyond
the hypothesis of interest. In order for his conclusion to
be relevant to the case, we must make the auxiliary assumptions:
\begin{itemize}
\item the probability of an incident during a night shift is the same
as during a day shift (but more people die during the night);
\item the probability of an incident during a shift does not
depend on the prevailing atmospheric conditions 
(but they may have an effect on respiratory problems);
\item the case-mix at the ward did not systematically change
over the time period concerned;
\item the occurrence of an incident in shift $n+1$ is
independent of the occurrence of an incident in shift $n$ (however, a
successful reanimation in shift $n$ may be followed by death in shift $n+1$);
\item in normal circumstances, all nurses have equal 
probability to witness incidents 
(on the contrary, as our own informal enquiries in hospitals 
have shown, terminally ill patients often die in the presence of a nurse 
with whom they feel `comfortable').
\end{itemize}
This is just a small sample of the auxiliary hypotheses
which are needed to make the rejection of 
Elffers' null-hypothesis relevant to the case at hand.
The main point is this:  only if the auxiliary
hypotheses used in setting up the model are realistic, can 
the occurrence of an improbable outcome be used to cast doubt on 
the null-hypothesis of interest. In the absence of such independent
verification of the auxiliary hypotheses, the occurrence of an 
improbable outcome might as well point to incorrect auxiliary
hypotheses. 

To put it a different way, Elffers' explicit model assumptions
show how he chose to formally interpret the question asked by
the court: could so many incidents occur during Lucia's shifts
by chance?  Our first three items listed above
suggest that the chance of an incident
might vary strongly over the shifts without there being any
difference between the nurses. Then if Lucia for whatever
reason tended to be assigned many more `heavy' shifts
than the other nurses, she will experience by chance 
much more than the average number of incidents. The assignment 
of nurses to shifts was certainly not done completely at random.
Notice that at the JKZ, Lucia had a much larger proportion of 
shifts than the other nurses. The nurses are not all the same 
in this respect.

Finally, as our last item shows, there may be completely
innocent reasons why the chance of an incident in a particular
shift might depend on the nurse who is on duty. 

\section{The court's interpretation of Elffers' numbers}

In its judgement of March 24, 2003, the court glossed Elffers' 
findings as follows (the numbering corresponds to the court's 
report; the emphasis -- \emph{italic script} -- is ours)\footnote{The 
original Dutch version can be found at www.rechtspraak.nl}:
\begin{quote} 7. In his report of May 29, 2002, Dr. H. Elffers concludes 
that the \textit{probability} that a nurse \textit{coincidentally} experiences 
as many incidents as the suspect is less than 1 over 342 million. 
\end{quote}
\begin{quote} 8. In his report of May 29, 2002, Dr. H. Elffers has further 
calculated the following component probabilities
\begin{enumerate}
\item[a.] The \textit{probability} that one out of 27 nurses would 
\textit{coincidentally} experience 8 incidents in 142 out of a total of 1029 
shifts \ldots is less than 1 over 300,000.
\item[b.] The \textit{probability} that the suspect has \textit{coincidentally} 
\ldots
\end{enumerate} 
\end{quote}
\begin{quote} 11. The court is of the opinion that the probabilistic 
calculations
given by Dr H. Elffers in his report of May 29, 2002, entail that it must be
considered \textit{extremely improbable} that the suspect  experienced all 
incidents mentioned in the indictment \textit{coincidentally}. These 
calculations \textit{consequently show} that it is \textit{highly probable} that
there is a \textit{connection between the presence of the suspect and the 
occurrence 
of an incident}. 
\end{quote}
We have cited these excerpts from the court's judgement because the italicised 
phrases should raise eyebrows among statisticians. The judgement of the court
is ambivalent, and it is unclear whether or not the court makes the
famous mistake known as the {\em prosecutor's fallacy}. Clearly, one {\em should} 
talk about the probability that something happened, under the assumption that
everything was totally random. The judgement of the court could however also be
interpreted as the probability that something accidentally happened. 
This is quite different, as is easily illustrated with the following 
formal translation into mathematical language. 

Writing $E$ for the observed event, and $H_0$ for the 
hypothesis of chance, Elffers calculated 
$P(E\mid H_0) < 342 \cdot 10^{-6}$, while the court seems to have 
concluded that $P(H_0\mid E) < 342 \cdot 10^{-6}$. Writing 
$$P(H_0\mid E) = \frac{P(E\mid H_0)\cdot P(H_0)}{P(E)},$$ 
we see that prior information about $P(H_0)$ and $P(E)$ 
would required to come to such a conclusion.

We would like to note that Elffers did not make this mistake himself, but during 
the testimony of two of the authors of this article (RM and MvL), 
the court of appeal certainly did. 

\section{Elffers' method revised}
\label{revise} 

There is an simple way to revise Elffers' method 
in a way that avoids the scale problems and the double use of data:
discard the data from JKZ, and just analyse the data from the
two wards at RKZ, combining the results in a statistically
correct fashion.

It was the concentration of incidents during Lucia's shifts at JKZ 
which suggested criminal activity with herself as suspect.
Elffers' analysis of those numbers informally
confirms that the concentration was surprising and 
justifies the further investigation that took
place. But the probability he reported to the court
for the JKZ is misleading, if not meaningless.
This does not mean that no evidence from 
the JKZ can be used in court; it just means that 
this particular data from the JKZ cannot be
used in a {\it statistical fashion}, at least,
not within the classical frequentist paradigm.
Data from the JKZ, for instance toxicological reports 
can be used in court in different ways.

Doing similar computations as Elffers, but now restricted to the RKZ and without
any correction, we obtain very different numbers. If we first take
the data of the two wards together, then we have a total number of $675$ 
shifts, Lucia having $61$ of them (note the correction of numbers). 
There were $19$ incidents, $6$ of which were during one
of Lucia's shifts. Under the same hypothesis as Elffers, a similar computation now 
leads to a probability of $0.0038$, which of course is much larger than the
number obtained by Elffers. In particular, Elffers himself used a significance level
of $0.001$, meaning that in this case the null-hypothesis should {\em not} be rejected, in 
sharp contrast to Elffers' conclusion.

However, one should make a distinction between the two wards,
which took rather different kinds of patients, and indeed
the rate of incidents in each seems quite different;
Lucia has proportionately more shifts in the ward where incidents
are more frequent. 
There are several ways of taking account of this. 
One can combine two separate $p$-values as in Section
\ref{mult}, or, alternatively, treat both wards independently 
with the hypergeometric method
as Elffers, and ask for the probability that the sum of 
two independent hypergeometric random variables 
(with their respective parameters) exceeds 6. 
A simple computation leads to the conclusion that this 
probability is equal to $0.022$, still bigger than the previously 
found $0.0038$.

It is clear that some of the aforementioned problems 
remain in this revised form of the method.  Nevertheless, 
we believe that the revised form is an improvement, 
since there is no double use of data, hence no need of
a post hoc correction without rationale. 
There are still subjective choices to be made
(how to combine the data from the two wards at RKZ)
but this is a matter of taste, not controversy. 
The revised analysis shows that the data from RKZ 
gives independent though rather weak confirmation 
that the rate of incidents was larger in Lucia's shifts 
than in those of other nurses.
This does not imply that her presence is the cause.
Without any information about the expected rate of
incidents, about how it might vary over different
kinds of shifts, and about how nurses are assigned to
different shifts, the data is rather inconclusive.

\section{A Bayesian approach to the problem}
\label{bayes}

During and after the trial, a public debate arose 
in The Netherlands about the way 
statistics was used in this case. Apart from Henk Elffers 
and two of the authors of this article, also Aart de Vos, 
an econometrician, entered the discussion. 
De Vos claimed that a Bayesian approach would solve 
scale problems and problems of post hoc data analysis; 
see \cite{vos1}-\cite{vos4}. In a national newspaper, 
he came to the conclusion that Lucia was {\em not} guilty 
with probability at least $10\%$, a number in sharp contrast 
with Elffers' outcomes. We summarise his method here, 
without going into details.

A Bayesian analysis works as follows. Let $E$ denote the evidence at hand,
$H_d$ the null-hypothesis (the hypothesis that L is innocent), and 
$H_p$ denote the alternative hypothesis (the hypothesis that L is guilty).

A straightforward application of Bayes' rule now gives
$$
\frac{P (H_p|E)}{P (H_d|E)}  = \frac{P
(E|H_p)}{P (E|H_d)}\cdot \frac{P (H_p)}{ P (H_d)}.
$$
In (other) words, 
$$
\mbox{posterior odds }= LR \cdot \mbox{ prior odds.}
$$
We interpret $P(H_d|E)$ as the probability of $H_d$ after evaluating 
the evidence $E$. The posterior odds are - at least in theory - nice 
to work with, because any new evidence ($E_{\mbox{new}}$) can be
implemented to give new posterior odds. For example, suppose we 
first had
$$
\mbox{``old'' posterior odds }=\frac{P
(E|H_p)}{P (E|H_d)}\cdot \frac{P (H_p)}{ P (H_d)},
$$ 
then, after this new evidence, we get new posterior odds:
\begin{eqnarray*}
\frac{P(H_p|E,E_{\mbox{\small new}})}{P(H_d|E,E_{\mbox{\small new}})}& = & 
\frac{P(E_{\mbox{\small new}}\cap
E|H_p)}{P
(E_{\mbox{\small new}}\cap E|H_d)}\cdot \frac{P (H_p)}{ P (H_d)}\\
& = & \frac{P(E_{\mbox{\small new}}|H_p,E)}{P(E_{\mbox{\small new}}|H_d,E)}\cdot
\mbox{``old'' posterior odds}.
\end{eqnarray*}

This is all nice in theory, but the questions that arise once you try to use
this in a law suit are obvious: can we make sense of $P(H_p)$ and $P(H_d)$? 
For what kind of evidence it is possible to compute
${P (E|H_p)}/{P (E|H_d)}$? And can we make sense of
${P(E_{\mbox{\small new}}|H_p,E)}/{P(E_{\mbox{\small new}}|H_d,E)}$? 
The latter question is particularly challenging, 
because it is difficult to see how the different 
pieces of evidence are related. 

In the case at hand, the following facts were brought up by De Vos 
as relevant evidence. After each piece of evidence we write between 
parentheses the likelihood ratio for that piece of evidence as used 
by De Vos.

\begin{enumerate}
\item $E_1$; the fact that the suspect never confessed ($\frac12$); 
\item $E_2$; the fact that two of the patients had certain toxic substances 
in their blood ($50$); 
\item $E_3$; the fact that 14 incidents occurred during Lucia's shifts ($7,000$);
\item $E_4$; the fact that suspect had written in her diary that 
`she had given in to her compulsion' ($5$).
\end{enumerate}

It seems obvious to us that these facts are hardly, 
if at all, expressible as numbers; the numbers of 
De Vos can hardly be justified. The prior probability 
$P(H_p)$ is taken to be $10^{-5}$, and then finally, 
De Vos assumes independence between the various facts, 
ending up with posterior odds equal to roughly $8.75$. 
This means that suspect is guilty with probability close
to $90\%$, certainly not enough to convict anybody.

\subsection{Discussion}

The numbers obtained by De Vos are in sharp contrast with Elffers' 
outcomes. However, it is clear from the analysis that his priors 
and likelihood ratios are very subjective. Any change in his 
priors would lead to very different answers. 

An advantage of the Bayesian approach is that there are no 
worries with post hoc corrections or scale problems: the priors 
should take care of these. Moreover, there are some constructive 
ideas in the modelling assumptions of De Vos; 
for instance, in order to arrive at a likelihood ratio 
for $E_3$, the number of incidents in Lucia's shifts, he 
proposes to take account of `normal' variation of incident rates 
between nurses, and he explains how he would estimate this
if relevant data were available (for now, he makes do with
a guess). On the other hand, he also has to come up with a 
probability for the precise number of incidents in Lucia's shifts 
if she is guilty! 

De Vos would like to see the Bayesian approach applied to
the case in its totality. The judge will base his verdict on
his posterior probability that the suspect is guilty.
This would require judges to give their priors in order 
to motivate their verdicts. It is unclear what the role of 
the defence would be in this situation: can they reasonably 
object to the judges' subjective priors? 

\section{An epidemiological approach}
\label{ait}

In \cite{lucy-aitken1} and \cite{lucy-aitken2}, 
Lucy and Aitken discuss a different way of
modelling cases like this, and we include a discussion of their method 
here. This method does not rely on conditioning on the number of 
incidents, but instead presumes availability of epidemiological data. 

The basic assumption of Lucy and Aitken is that the probability 
distribution of the number $X$ of incidents witnessed by a certain nurse, 
is given by a Poisson distribution, hence
$$
P(X=k)=e^{-\mu r}\frac{(\mu r)^k}{k!},
$$
where $r$ is the number of shifts of the nurse, and $\mu>0$ is a parameter 
representing the intensity of incidents. 

The usual argument for the Poisson distribution in this kind of
situation is that it follows from the following assumptions:
the numbers of incidents in different time intervals are 
independent of one another, with constant expected rate; 
several incidents can not occur at the same time. 
Since the chance of several incidents in one shift is rather
small, these assumptions are very close to those made by Elffers.
Indeed, the binomial distribution with small `success' probability
$p$ is very close to a Poisson distribution.
(One can therefore make the same objections to this model as to
Elffers': is the incident rate constant, are incidents at different
times independent of one another?)

The hypothesis of chance could now be formulated as saying that 
every nurse, including the suspect, has the {\em same} 
intensity parameter $\mu$. (Aart de Vos would allow every nurse to
have a \emph{different} intensity; the incident intensities of
innocent nurses being drawn from some probability distribution.)
The hypothesis $H_p$ of the prosecutor can have several forms. 
One possibility is that incidents in Lucia's shifts also follow
a Poisson distribution, but with a different intensity.
Then the prosecutor's hypothesis might be $H_p: \mu_L > \mu$, 
where $\mu_L$ is the parameter corresponding
to the suspect, and $\mu$ is the parameter corresponding 
to all other nurses, neither being specified. 
How to proceed, depends on whether or not $\mu$ and/or
$\mu_L$ are known or unknown quantities.

\subsection{Likelihood ratios}

One possible approach is to compute likelihood ratios for 
$H_p$ against $H_d$. Consider a situation with $I$ nurses, 
and let $k_i$ be the number of incidents witnessed by nurse 
$i$, $i=1,\ldots, I$. Denote by $r_i$ the
number of shifts of nurse $i$, and let $E$ be the event that
nurse $i$ witnessed $k_i$ incidents, for $i=1,\ldots, I$.  
This leads to
$$
P(E|H_d)=\prod_{i=1}^I e^{-\mu r_i} \frac{(\mu r_i)^{k_i}}{k_i!},
$$
and, assuming that the suspect is nurse $j$, to
$$
P(E|H_p)=\frac{e^{-\mu_L r_j}(\mu_L r_j)^{k_j}}{k_j!}
\prod^I_{i=1, i\neq j}e^{-\mu r_i} \frac{(\mu r_i)^{k_i}}{k_i!}.
$$
A simple computation that shows that the likelihood ratio becomes
\begin{equation}
\mbox{LR}=\frac{P(E|H_p)}{P(E|H_d)}=e^{\mu r_j-\mu_L r_j}(
\frac{\mu_L r_j}{\mu r_j} )^{k_j}.
\end{equation}

In order to evaluate the outcome of any computation with this likelihood ratio, 
we may use the following scale for describing the height of a likelihood ratio, 
see \cite{evett}: 

\medskip\noindent
\begin{tabular}{|c|c|}
\hline
 & evidence is\\
\hline
$LR=1$& equally likely under $H_p$ as under $H_d$\\
\hline
$1<LR<100$& slightly more likely under $H_p$ than under $H_d$\\
\hline
$100\leq LR<1000$& more likely under $H_p$ than under $H_d$\\
\hline
$1000\leq LR<10,000$& much more likely under $H_p$ than under
$H_d$\\
\hline
$LR>10,000$& very much more likely under $H_p$ than under
$H_d$\\
\hline
\end{tabular}

\medskip\noindent

As was noted by Meester and Sjerps in \cite{ms1} and \cite{ms2}, one should be
careful when using a table like this if the hypotheses were suggested by 
the data. In that case they only become meaningful in combination with
prior probabilities for the hypotheses considered. For this reason we concentrate 
on the RKZ. However we still have a problem with data dependent hypotheses, since
we need to specify the intensities $\mu$ and $\mu_L$ in order to compute
the likelihood ratio.

In the following computations, for simplicity we take
the data of the two wards at the RKZ together. Above we have argued that 
we should allow different incident rates between different wards; 
in that case the numbers would come out even better for the suspect.

\medskip\noindent
{\bf I:} Without further data, a reasonable assumption for the prosecutor 
is to estimate $\mu$ using the incidents during shifts of all nurses
apart from the suspect. 
$$
\mu=\frac{13}{614}.
$$
Lucy and Aitken proceed by choosing $\mu_L$ in such a way that the expected number
of incidents witnessed by the suspect is precisely $k_j$, that is, $\mu_L r_j=k_j$ 
hence
$$
\mu_L=\frac{6}{61}.
$$
These assumptions lead to a likelihood ratio of $90.7$, and this is in the range
where the evidence is only slightly more likely under $H_d$. 

\medskip\noindent
{\bf II:} The defence might prefer to estimate $\mu$ based on all incidents, 
we would then get 
$$
\mu=\frac{19}{675},
$$
and this leads to a likelihood ratio of about $25$ (keeping $\mu_L$ as above).

If we would apply this method to the JKZ data, prosecution and defence 
would disagree strongly on how to estimate $\mu$. From the point of view of
the defence, the prosecution's estimate is biased downwards, grossly.
The precise reason we are analysing this data, is because we observed 
a coincidental concentration of incidents in the shifts of one nurse; 
we then take the other shifts, with coincidentally few incidents, 
on which to base our estimate!

\subsection{Relation to Elffers' approach}

As we noted above, the assumptions needed to justify the Poisson 
model are essentially the same as those \emph{initially} taken by Elffers.
Starting from the model of Lucy and Aitken, conditioning on the
observed fact that there was never more than one incident in a shift,
and then conditioning on the total number of incidents, 
we arrive, under the null-hypothesis of chance, at 
Elffers' hypergeometric distribution. Now, one can 
also arrive at the hypergeometric distribution from different 
modelling assumptions; for instance: the chance of an 
incident in a shift may vary arbitrarily over shifts, 
but nurses are assigned to shifts completely at random.
Thus Lucy and Aitken's analysis is more restrictive that
that of Elffers. In particular, the Poisson model
suffers from all the problems brought up in Section \ref{qd}; 
also in this model, we must assume there is no difference
between day and night shifts, and no variation in case-mix 
over time, and so on.

If the `normal' intensity were known to be equal to $\mu$, 
then using the property of sufficiency we see that
the data of the other nurses is irrelevant and we should simply
investigate whether Lucia's number of incidents is large compared to
the number expected from a Poisson distribution with expectation
$\mu$ times the number of shifts of Lucia.

If the normal intensity is unknown but we have data from `normal'
working operations (e.g.\ the other nurses in the same ward) then
again sufficiency shows that we should base our inference on the
total numbers of incidents of Lucia on the one hand, and of the 
others on the other hand. These numbers will
be Poisson distributed with means $\mu_L r_L$ and $\mu r$ 
respectively, where $r_L$ and $r$ are the number of shifts 
of Lucia and the others, respectively.

This is a classical statistical hypothesis testing problem.
If both means are large, one would use a generalised likelihood
ratio procedure based on comparing maximised log likelihoods under
the null-hypothesis: $\mu_L=\mu$, both parameters unknown.
However, many statisticians would prefer to use an exact test 
based on the fact that conditional on the grand total of 
incidents $N$, those of Lucia are binomially distributed 
with parameters $N$ and 
$$
p=\frac{\mu_L r_L}{\mu_L r_L + \mu r}.
$$
Under our null-hypothesis $\mu=\mu_L$, $p$ is known, and we have 
a classical hypothesis testing problem based on one observation from a 
binomial distribution. If we have no more data than that reproduced 
in the paper then (per ward) the analysis is almost the same as Elffers'.

However, if there truly were more data available, e.g.\ numbers of incidents 
and shifts in some adjacent time periods in the same ward, then this has
the effect of adding to the total number of shifts and adding to the
total number of incidents during the shifts of the others. If we had 
much of such data and if the incident rate during that time was close to 
what we observed during Lucia's shifts, then the data would
become more and more favourable to Lucia. 
Hence, {\it if} we had had information on the `normal' rate of
incidents, we would have used it, and the conclusion {\it could} 
have been very different. 

\subsection{Discussion}

The drawbacks of the conditional approach become quite
apparent here. In the previous subsection it became 
clear that collecting extra data could change the impact
of the existing data dramatically, and it would seem 
the duty of an expert witness to point this out.

However, if reliable data on incident rates cannot be found, 
this approach leads to essentially the same analysis and
conclusions as the corrected Elffers method, 
see Section \ref{revise}. Also the likelihood 
ratio approach, applied to just the RKZ data, 
leads to much the same conclusion again.

\section{Relative risk}
\label{rr}

In \cite{lucy-aitken1} and \cite{lucy-aitken2}, 
Lucy and Aitken define the term \emph{relative risk} as
follows: the relative risk $R_j$ of a nurse $j$ 
is the fraction of her shifts during which an incident
took place, divided by the fraction of the remaining shifts 
during which an incident took place. More formally,
$$
R_j=\frac{k_j/r_j}{\sum_{i \neq j} k_i /\sum_{i \neq j} r_i}.
$$ 
For example, the relative risk of Lucia for the RKZ
for the two wards together is equal to 
$$
\frac{\frac{6}{61}}{\frac{13}{614}}\approx 4.65.
$$

The fact that Lucia had the highest relative risk is clearly not enough
to warrant any investigation; some nurse must have the highest relative risk.
The more important question is how high a relative risk should be in order
to be suspicious. 

The distribution of the highest relative risk depends on many 
variables, like the number of nurses, the way the shifts are spread 
among the nurses, the number of shifts, and of course on 
the modelling assumptions concerning the occurrence of incidents. 
In this section we again concentrate on the model of Lucy and Aitken 
of the previous section.

The numbers in the definition of the relative risk only
depend on the considered time span, comparing the amount of
incidents (s)he witnessed to the amount of incidents the 
other nurses witnessed. It is now useful to do some numerical 
simulations to obtain some idea about the distribution of 
the highest relative risk.

\subsection{Simulating relative risk}

We have no data concerning the number of shifts of each other nurse, 
apart from Lucia; and at the RKZ we do not even know how many other
nurses there were. Therefore, we have simulated a situation where all nurses 
worked the same number of shifts (actually, this should lead to less
variability in relative risk, since if some nurses work few shifts, their
relative risks can more easily be extremely small or large). 
At the RKZ we have a total of 675 shifts of which Lucia worked 61 
(note the remark after the table of data in Section 2). 
Therefore we simulated a situation in which 11 nurses
all had 61 shifts. Hence $r_i=r$ for all $i$ and $I=n/r$. 
This leads to 
$$
R_j= \frac{k_j}{\sum^I_{i=1}k_i-k_j}(I-1).
$$
We are interested in the nurse with the highest relative risk for each group
of $I$ nurses. Since all nurses work the same number of shifts, this is simply
the nurse with the most incidents.

We have run $1000$ simulations in the case of Lucia for the data of the RKZ, 
first for both wards together, then for each ward separately. 
The values for $\mu$ in the first column are based on
the frequency of incidents of all other nurses in the RKZ; 
the values of $\mu$ in the second column are based on the 
overall frequency of incidents, including Lucia 
(these choices are the same as in Section 7.1).

\medskip\noindent
\begin{tabular}{|l|l|l|}
\hline
 whole RKZ& $\mu=\frac{13}{614}$ &
$\mu=\frac{19}{675}$ \\
Lucia's $p$-value & 0.121 & 0.042  \\
\hline
RKZ-41 & $\mu=\frac{4}{333}$ & $\mu=\frac{5}{336}$ \\
Lucia's $p$-value & 0.787 & 0.681 \\
\hline
RKZ42 & $\mu=\frac{9}{281}$ & $\mu=\frac{14}{339}$ \\
Lucia's $p$-value & 0.383 & 0.286 \\
\hline
\end{tabular}\\[1mm]

\subsection{Discussion}

For $\mu=\frac{13}{614}$, L's relative risk of approximately $4.65$ 
lies between the $879$th and $880$th of the $1000$ highest relative
risks. In other words, it is a high relative risk, but not extremely high.
For $\mu=\frac{19}{675}$, L's relative risk lies between the $958$th and the
$959$th highest relative risks. If we would take $\mu$ even higher, 
L's relative risk would have a smaller $p$-value.

 From this, we may conclude that if data on the number of incidents outside
the time span L worked at the RKZ would indicate $\mu$ to be
large, L's relative risk would be extremely high and this could be used as
evidence against her in court. This seems strange, 
since in the likelihood ratio approach of the previous section, 
a larger $\mu$ implied a {\em lower} likelihood ratio, 
which is in favour of the defendant (if the Poisson model 
is correct and $\mu$ is known, then only the number of 
incidents in Lucia's shifts is relevant for investigating
whether her incidents have the same intensity).

The fact that a large $\mu$ does not work
in favour of the defendant in the relative risk approach, is because if 
$\mu$ really is very large, then we do not expect much spread in the the
relative risks of the nurses. If the total number of incidents is
coincidentally very small, then the relative risks will be widely spread.
So under the hypothesis that all the nurses are the same, whenever
the total number of incidents is much smaller than expected, 
the largest relative risk is likely to be extremely large.

\section{Conclusion}
It is not easy to draw a clear cut conclusion from all this. 
Elffers' analysis of the JKZ data perhaps confirms that something
surprising has happened there, just as you would be surprised 
if someone in your street won the state lottery.
Indeed, had this computation led to the conclusion 
that the concentration of incidents in Lucia's shifts was not so
surprising, then there would not have been a case against 
Lucia at all. 

So, these numbers did raise interest and suspicion, and 
there should have been reflection on what to do next. 
If Elffers used his model correctly, that is, 
combining the data from different wards in a statistically
justifiable way; moreover, without double use of data 
(hence without the need for arbitrary  post hoc correction), 
then the resulting numbers would have been very different. 
In fact, the outcome would not have led to the rejection of
the null-hypothesis of chance, at significance level of $0.001$ 
(Elffers' own choice),  although the $p$-value of $0.022$
(see Section \ref{revise}) still leaves one uneasy. 

Following the epidemiological approach does not lead to a different 
conclusion. The likelihood ratios of $90.7$ and $25$ reported in Section 
\ref{ait} would in itself not lead to conviction, but are again uneasily 
high. Similar remarks apply to the relative risks in Section \ref{rr}. 

De Vos' Bayesian approach does combine all the data. He starts with a
rather small prior probability that Lucia is a murderer. He allows
natural variation between incident rates between innocent nurses,
making Lucia's number of incidents somewhat less surprising
(a likelihood ratio of `only' $7,000$ for the JKZ data). 
Together with a small prior probability for Lucia to be a
murderer, he arrives at a chance of $10\%$ that she is innocent.
However on the way he has to conjure up one number after another
out of thin air.

In contrast to this, the weakness of Elffers' approach can be seen as 
its strength.  Once we have the model, there are no further parameters 
to be worried about or which could lead to disagreement between
prosecution and defence. Convincing rejection of the null-hypothesis 
would mean that the association between incidents and Lucia's shifts 
is not a coincidence. However, correlation does not imply causation,
and alternative explanations of the correlation need to be disqualified
before it could be seen as evidence for the case of the prosecution.
(The conclusion of the court was that there was a `connection' between her
presence, and incidents; the word \emph{connection} is oversuggestive
of causality).

The Poisson model of Lucy and Aitken suffers from the fact
that any conclusion by either party can be questioned 
by the other on the basis of the choice of the parameter $\mu$. 
If one of the parties can raise reasonable doubts about the 
validity or reasonableness of the parameter choice, 
then the numbers arising from that model can be 
questioned as well. 

On the other hand, the analysis carried out in Section \ref{ait} 
shows that Elffers' choice to condition on the total number of incidents
is tantamount to ignoring what could be a very relevant piece of
information.

The more sophisticated a model becomes, the more possibilities for
criticising it one has. This becomes abundantly clear in the 
Bayesian approach of De Vos in Section \ref{bayes}. 
De Vos tries to incorporate everything into his mathematical model. 
To us, this seems impossible, and the result of the computations
of De Vos do not mean much.

Can statistics play an important role in a case 
like this? As we have seen there is no one
correct way to analyse the available data. Every analysis involves
subjective choices. The more sophisticated the analysis, the more
subjective elements it seems to contain, and hence the more
controversial are its conclusions.
In fact, perhaps the only uncontroversial number in this paper 
are the numbers in Section \ref{revise}. These numbers do not 
suffer from double use of data or scaling problems, nor do
they involve any parameter choice. On the other hand, the evidence 
they give is weak, in several respects (weak evidence of correlation,
not of causation).

Perhaps these numbers are -- after all -- the only possible contribution 
of statistics to the present case. This statement might be surprising, 
given the data in the table at the beginning of the paper. 
But it is one thing to say that a number is relevant,
it is quite another thing to work out a reasonable way to use it. 
Not all numbers can or should be used in a statistical fashion.
  
On June 18, 2004, the court of appeal in The Hague again
found Lucia de B. guilty and sentenced her to life imprisonment plus detention 
in a psychiatric hospital in case she would ever be pardoned. This time the
judgement made no mention at all of statistical arguments; other evidence which
had played a secondary role during the first trial now assumed primary
importance. Hence for several reasons this was a Pyrrhic victory at most
(at least for the authors). If the court of appeal had explicitly repudiated
the form of statistical argument employed by the public prosecutor and the
first court, future cases would have been able to use this jurisprudence.

However, incorporating the statistical argument in the second judgement would
have required the court of appeal to take an explicit stand on all the 
issues raised above. In fact, careful writers on the foundations of 
statistics have pointed out that  evaluating a statistical conclusion involves 
even more:
\begin{quote} In applying a particular technique in a practical problem,
it is vital to understand the philosophical and conceptual attitudes
from which it derives if we are to be able to interpret (and appreciate 
the limitations of) any conclusions we draw. 
(\cite{barn}, page 332)\end{quote}
Evidently the court of appeal was not willing to dig this deep;
but the quote as well as the case of Lucia de B. may serve as 
a reminder to lawyers and judges that the interpretation of statistical 
arguments is by no means immune to disputation.

\end{document}